\documentclass[11pt]{article}
\usepackage{amsmath}
\usepackage{amssymb}
\def\C{\mathbb{C}}
\def\R{\mathbb{R}}
\def\N{\mathbb{N}}
\newtheorem{theorem}{\hspace*{\parindent}Theorem}

\newtheorem{corollary}{\hspace*{\parindent}Corollary}

\title{Representations and inequalities for generalized hypergeometric functions}
\author{Dmitrii Karp\footnote{Far Eastern Federal University, Vladivostok, Russia,
e-mail:\,\emph{dimkrp@gmail.com}}}
\date{}
\begin{document}
\maketitle

\begin{center}
\parbox{12cm}{
\small\textbf{Abstract.}
We find an integral representation for the generalized hypergeometric function unifying known representations via generalized Stieltjes, Laplace and cosine Fourier transforms. Using positivity conditions for the weight in this representation we establish  various new facts regarding generalized hypergeometric functions, including complete monotonicity, log-convexity in upper parameters, monotonicity of ratios and new proofs of Luke's bounds.  Besides, we derive two-sided inequalities for the Bessel type hypergeometric functions by using their series representations.}
\end{center}

\bigskip


\bigskip


\bigskip

\section{Introduction}

We will adopt standard notation $\R$, $\C$ and $\N$  for the real, complex and natural numbers, respectively.  $\N_0$ will denote $\N\cup\{0\}$. In our previous works \cite{KPJMAA,KSJAT} we obtained some representations, inequalities, monotonicity and other properties for the Gauss type generalized hypergeometric function $_{q+1}F_q$ which is equal to $p=q+1$ case of the function \cite{NISTCh16,Kiryakova}
\begin{equation}\label{eq:pFqdefined}
{_{p}F_q}\left(\left.\!\!\begin{array}{c} A\\B\end{array}\right|z\!\right)={_{p}F_q}\left(A;B;z\right):=
\sum\limits_{n=0}^{\infty}\frac{(a_1)_n(a_2)_n\cdots(a_{p})_n}{(b_1)_n(b_2)_n\cdots(b_q)_nn!}z^n,
\end{equation}
where $A=(a_1,a_2,\ldots,a_p)$ and $B=(b_1,b_2,\ldots,b_q)$, $b_j\notin-\N_0$, are parameter vectors, $(a)_n$ denotes the rising factorial, defined by $(a)_0=1$, $(a)_n=a(a+1)\cdots(a+n-1)$, $n\geq{1}$. The series in (\ref{eq:pFqdefined}) converges in the entire complex $z$-plane if $p\leq{q}$ and inside the unit disk if $p=q+1$. In the latter case the sum can be extended to a function holomorphic in the cut plane  $\C\!\setminus\![1,\infty)$.  The main tool employed in \cite{KPJMAA,KSJAT} to investigate the function $_{q+1}F_q$ is the generalized Stieltjes transform (see (\ref{eq:FStieltjes}) below) of a measure with density expressed by the $G$-function of Meijer, cf. \cite[Theorem~2]{KPJMAA}. Such representation appeared earlier in \cite[Theorem~4.2.11]{Kiryakova}.  We contributed more relaxed conditions on parameters and studied nonnegativity of the representing measure.  This lead to monotonicity of the ratios, two-sided bounds, mapping properties and other results for the Gauss type hypergeometric functions $_{q+1}F_q$.

Another line of research pursued in \cite{KKITSF,KKJMAA,KSJMAA} hinges on the series representation (\ref{eq:pFqdefined}) and yields, among other things, a number of properties of the Kummer type hypergeometric functions ${_qF_q}$, including logarithmic concavity or convexity in parameters, inequalities for logarithmic derivatives and bounds for the Tur\'{a}nians.  In this note we introduce an integral representation for the general hypergeometric function ${_pF_q}$, which includes, as particular cases, the representations by the generalized Stieltjes, Laplace and cosine Fourier transforms. Starting with this representation we will  obtain new properties of the the Gauss type functions ${_{q+1}F_q}$, the Kummer type functions $_qF_q$ and the Bessel type functions $_{q-1}F_q$, including conditions for complete monotonicity, monotonicity of ratios and log-convexity in upper parameters. Moreover, we furnish new proofs for Luke's inequalities from \cite{Luke}, allowing their extension to a wider parameter range. Finally,  we discover new bounds for the Bessel type hypergeometric functions $_pF_q$ with $p<q$ of positive argument.

\section{Representations for ${_pF_q}$ and their consequences}
Suppose $0\leq{m}\leq{q}$, $0\leq{n}\leq{p}$ are integers and $A\in\C^{p}$,  $B\in\C^{q}$ are such that  $a_i-b_j-1\notin\N_0$ for all $i=1,\ldots,n$ and $j=1,\ldots,m$.  We will heavily use Meijer's $G$-function \cite[Section~16.17]{NISTCh16} defined by the contour integral
\begin{equation}\label{eq:G-defined}
G^{m,n}_{p,q}\!\left(\!z~\vline\begin{array}{l}A
\\B\end{array}\!\!\right)\!\!:=
\frac{1}{2\pi{i}}
\int\limits_{\mathcal{L}}\!\!\frac{\Gamma(b_1\!+\!s)\cdots\Gamma(b_m\!+\!s)\Gamma(1-a_1\!-\!s)\cdots\Gamma(1-a_n\!-\!s)z^{-s}}
{\Gamma(a_{n+1}\!+\!s)\cdots\Gamma(a_p\!+\!s)\Gamma(1-b_{m+1}\!-\!s)\cdots\Gamma(1-b_{q}\!-\!s)}ds.
\end{equation}
The contour $\mathcal{L}$ begins and ends at infinity and separates the poles of the integrand of the form $-b_j-k$, $k\in\N_0$, leaving them on the left, from the poles of the form  $-a_j+k+1$, $k\in\N_0$, leaving them on the right. Under the above conditions such contour always exists and can be chosen to make the integral in (\ref{eq:G-defined}) convergent.  More details regarding the choice of $\mathcal{L}$ and conditions for convergence in  (\ref{eq:G-defined}) can be found in \cite{NISTCh16}, \cite[Chapters~1 and 2]{KilSaig} and \cite[Chapter~8]{PBM3}.

We will abbreviate $\prod_{i=1}^{p}\Gamma(a_i)$ to $\Gamma(A)$ and $\prod_{i=1}^{p}(a_i)_n$ to $(A)_n$ throughout the paper.  Expressions like $A+\alpha$, where $\alpha\in\C$, and $\Re(A)>0$ will be understood element-wise.  The key role in the investigations carried out in \cite{KPJMAA,KSJAT} is played by the generalized Stieltjes transform representation
\begin{equation}\label{eq:FStieltjes}
{_{q+1}F_q}\!\left(\left.\!\!\begin{array}{c}\sigma,A\\
B\end{array}\right|-z\!\right)
=\frac{\Gamma(B)}{\Gamma(A)}\!\int\limits_{0}^{1}(1+zt)^{-\sigma}G^{q,0}_{q,q}\left(t\left|\begin{array}{l}\!\!B\!\!\\\!\!A\!\!\end{array}\right.\right)\frac{dt}{t},
\end{equation}
which is easy to prove by termwise integration.  Note that both the generalized Stieltjes kernel $(1+zt)^{-\sigma}={_1F_0}(\sigma;-;-zt)$ and  the Laplace kernel $e^{-zt}={_0F_0}(-;-;-zt)$ are particular cases of a more general hypergeometric kernel.  This simple observation leads to the following theorem.

\begin{theorem}\label{th:pFqgeneral}
Suppose $p_1\ge0$, $p_2\ge1$, $q_1,q_2\ge0$, $p_2\ge{q_2}$, $p=p_1+p_2$, $q=q_1+q_2$, $p\le{q+1}$ are integers
\emph{(}these conditions imply that $p_1\leq q_1+1$\emph{)}. Write $A_1=(a_1,\ldots,a_{p_1})$, $A_2=(a_{p_1+1},\ldots,a_{p})$, $B_1=(b_1,\ldots,b_{q_1})$, $B_2=(b_{q_1+1},\ldots,b_{q})$ for complex parameter vectors satisfying $\Re(A_2)>0$. Then
\begin{equation}\label{eq:pFqgeneral}
{_{p}F_q}(A_1,A_2;B_1,B_2;-z)=\frac{\Gamma(B_2)}{\Gamma(A_2)}
\int\limits_{0}^{\infty}{_{p_1}\!F_{q_1}}(A_1;B_1;-zt)G^{p_2,0}_{q_2,p_2}\left(t\left|\begin{array}{l}\!\!B_2\!\!\\\!\!A_2\!\!\end{array}\right.\right)\frac{dt}{t}.
\end{equation}
This formula is valid for $z\in\C$ if $p_1\leq q_1$ or $z\in\C\!\setminus\!(-\infty,-1]$ if $p_1=q_1+1$;
if $p_2=q_2$ additional assumption $\Re(\psi_2)>0$, where $\psi_2=\sum_{i=p_1+1}^{p}(b_i-a_i)$, has to be adopted \emph{(}in this case the $G$-function in \emph{(\ref{eq:pFqgeneral})} vanishes for $t>1$\emph{)}. If  $p_2=q_2$  and $\psi_2=0$, then
\begin{equation}\label{eq:pFqpsi2-0}
{_{p}F_q}(A_1,A_2;B_1,B_2;-z)=\frac{\Gamma(B_2)}{\Gamma(A_2)}
\biggl\{{_{p_1}\!F_{q_1}}(A_1;B_1;-z)+\int\limits_{0}^{1}{_{p_1}\!F_{q_1}}(A_1;B_1;-zt)G^{p_2,0}_{q_2,p_2}\left(t\left|\begin{array}{l}\!\!B_2\!\!\\\!\!A_2\!\!\end{array}\right.\right)\frac{dt}{t}\biggr\},
\end{equation}
where $z\in\C$ if $p_1\leq q_1$ or $z\in\C\!\setminus\!(-\infty,-1]$ if $p_1=q_1+1$.
\end{theorem}
\textbf{Proof.} Once the correctness of termwise integration has been justified, it suffices to write the kernel
${_{p_1}\!F_{q_1}}$ as the series (\ref{eq:pFqdefined}) and  integrate term by term to establish formula (\ref{eq:pFqgeneral}).  To demonstrate convergence of the integral  in (\ref{eq:pFqgeneral}) and justify the exchange of summation and integration  we resort to the asymptotic relation
\begin{equation}\label{eq:G-asymp-zero}
G^{p_2,0}_{q_2,p_2}\!\left(x~\vline\begin{array}{l}\!B\!\\\!A\!\end{array}\right)=
O\left(x^{a}\ln^{m-1}(x)\right)~~\text{as}~~x\to{0},
\end{equation}
where $a=\min(\Re(a_1),\ldots,\Re(a_p))$, and the minimum is taken over those $a_i$ for which $a_i-b_j\notin\N_0$ for all $j=1,\ldots,q_2$. Positive integer $m$ is the maximal multiplicity among the numbers $a_i$ for which the minimum is attained. This formula follows from \cite[Corollary~1.12.1]{KilSaig} or \cite[formula~(11)]{Karp}.  It proves convergence in (\ref{eq:pFqgeneral}) around zero. Near infinity for $p_2>q_2$ we have
\begin{equation}\label{G-asymp-inf}
G^{p_2,0}_{q_2,p_2}\left(x~\vline\begin{array}{l}B\\A\end{array}\right)=
\frac{(2\pi)^{\frac{1}{2}(\mu-1)}}{\sqrt{\mu}}x^{(1-\alpha)/\mu}e^{-\mu
x^{1/\mu}}\left[1+O(x^{-1/\mu})\right]~~\text{as}~~x\to{\infty},
\end{equation}
where $\mu=p_2-q_2$, $\alpha=\sum_{i=q_1+1}^{q}b_i-\sum_{i=p_1+1}^{p}a_i+\frac{1}{2}(p_2-q_2+1)$.
This formula is a particular case of the formula on page 289 in \cite{Braaksma} which is implied by formula (7.8) of the same paper.  If $p_2=q_2$ and $\Re(\psi_2)>0$ then
$$
G^{p_2,0}_{q_2,p_2}\left(x\left|\begin{array}{l}\!\!B_2\!\!\\\!\!A_2\!\!\end{array}\right.\right)=O((1-x)^{\Re(\psi_2)-1})~~\text{as}~x\uparrow1
$$
according to \cite[8.2.59]{PBM3} and
$$
G^{p_2,0}_{q_2,p_2}\left(x\left|\begin{array}{l}\!\!B_2\!\!\\\!\!A_2\!\!\end{array}\right.\right)=0~~\text{for}~x>1,
$$
according to \cite[Lemma~1]{KPJMAA} (see also proof of Theorem~\ref{th:G-positive} below).  This shows the convergence in (\ref{eq:pFqgeneral}) around unity for $p_2=q_2$. Finally, formula (\ref{eq:pFqpsi2-0}) follows from \cite[Theorem~1]{KL}.$\hfill$$\square$

\noindent\textbf{Remark}. Condition $p_2\ge{q_2}$ is necessary in the above theorem since for при $p_2<{q_2}$
$$
G^{p_2,0}_{q_2,p_2}\left(x\left|\begin{array}{l}\!\!B_2\!\!\\\!\!A_2\!\!\end{array}\right.\right)=0~\text{for all}~x\in\R.
$$
This condition shows that for $p<q$ the function ${_{p}F_q}$ cannot be represented by the Laplace or generalized Stieltjes transform.   The most ''extreme'' representation we can get in this case is:
$$
{_{p}F_q}(A;B;-z)=\frac{\Gamma(B_2)}{\Gamma(A)}
\int\limits_{0}^{1}{_0F_{q-p}}(-;B_1;-zt)
G^{p,0}_{p,p}\left(t\left|\begin{array}{l}\!\!B_2\!\!\\\!\!A\!\!\end{array}\right.\right)\frac{dt}{t},
$$
where the kernel ${_0F_{m}}$ is essentially the Bessel function if $m=1$ or the so called hyper-Bessel function if $m>1$ (see \cite{Kiryakova}).  Besides, this kernel cannot be represented by Theorem~\ref{th:pFqgeneral} due to condition $p_2\ge1$.
 It is sometimes desirable, however, to have a representation with a kernel independent of the parameters of the function being represented. This can be easily achieved by introducing artificial parameters $\alpha_j>0$ to get
\begin{equation}\label{eq:pFqsmallp}
{_{p}F_q}(A;B;-z)=\frac{\Gamma(B)}{\Gamma(A)\prod_{i=1}^{q-p}\Gamma(\alpha_i)}
\int\limits_{0}^{1}{_0F_{q-p}}(-;\alpha_1,\ldots,\alpha_{q-p};-zt)
G^{q,0}_{q,q}\left(t\left|\begin{array}{l}\!\!B\!\!\\\!\!A,\alpha_1,\ldots,\alpha_{q-p}\!\!\end{array}\right.\right)\frac{dt}{t}.
\end{equation}
We need to require $\sum{b_i}>\sum{a_i}+\sum{\alpha_i}$ for convergence of the above integral.  In particular, choosing $\alpha_i=i/(q-p+1)$, we obtain the kernel in terms of the so called generalized cosine,
$$
\cos_n(z)=\sum\limits_{j=0}^{\infty}\frac{(-1)^jz^{nj}}{(nj)!}={_0F_{n-1}}(-;1/n,2/n,\ldots,(n-1)/n;-(z/n)^n).
$$
The representation with such kernel has been first suggested by Kiryakova in \cite{Kiryakova}. An important particular case $p=q-1$ leads to standard cosine kernel as indicated in the corollary below. Before stating it let us define the the parametric excess by
\begin{equation}\label{eq:psi-defined}
\psi=\sum\limits_{k=1}^{q}b_k-\sum\limits_{k=1}^{p}a_k.
\end{equation}

\begin{corollary}
Suppose $\Re(A)>0$ element-wise. Then
\begin{equation}\label{eq:q+1FqLaplace}
{_{q+1}F_q}\!\left(\left.\!\!\begin{array}{c}A\\
B\end{array}\right|-z\!\right)=\frac{\Gamma(B)}{\Gamma(A)}\int\limits_{0}^{\infty}e^{-zt}G^{q+1,0}_{q,q+1}\left(t\left|\begin{array}{l}\!\!B\!\!\\\!\!A\!\!\end{array}\right.\right)\frac{dt}{t}.
\end{equation}
If also $\Re(\psi)>0$, then
\begin{equation}\label{eq:qFqLaplace}
{_{q}F_q}\!\left(\left.\!\!\begin{array}{c}A\\
B\end{array}\right|-z\!\right)
=\frac{\Gamma(B)}{\Gamma(A)}\int\limits_{0}^{1}e^{-zt}G^{q,0}_{q,q}\left(t\left|\begin{array}{l}\!\!B\!\!\\\!\!A\!\!\end{array}\right.\right)\frac{dt}{t}.
\end{equation}
If $\Re(\psi)>1/2$, then
\begin{equation}\label{eq:q-1Fqcosine}
{_{q-1}F_q}\left(\left.\!\!\begin{array}{c}A\\B\end{array}\right|-z\!\right)
=\frac{\Gamma(B)}{\sqrt{\pi}\Gamma(A)}\int\limits_{0}^{1}\cos(2\sqrt{zt})G^{q,0}_{q,q}\left(t\left|\begin{array}{l}\!\!B\!\!\\\!\!A,1/2\!\!\end{array}\right.\right)\frac{dt}{t}.
\end{equation}
If $\psi=0$, then \emph{(}\ref{eq:qFqLaplace}\emph{)} takes the form
\begin{equation*}
{_{q}F_q}\!\left(\left.\!\!\begin{array}{c}A\\
B\end{array}\right|-z\!\right)
=\frac{\Gamma(B)}{\Gamma(A)}\biggl\{e^{-z}+
\int\limits_{0}^{1}e^{-zt}G^{q,0}_{q,q}\left(t\left|\begin{array}{l}\!\!B\!\!\\\!\!A\!\!\end{array}\right.\right)\frac{dt}{t}\biggr\}.
\end{equation*}
If $\psi=1/2$, then \emph{(}\ref{eq:q-1Fqcosine}\emph{)} takes the form
\begin{equation*}
{_{q-1}F_q}\left(\left.\!\!\begin{array}{c}A\\B\end{array}\right|-z\!\right)
=\frac{\Gamma(B)}{\sqrt{\pi}\Gamma(A)}\biggl\{\cos(2\sqrt{z})
+\int\limits_{0}^{1}\cos(2\sqrt{zt})G^{q,0}_{q,q}\left(t\left|\begin{array}{l}\!\!B\!\!\\\!\!A,1/2\!\!\end{array}\right.\right)\frac{dt}{t}\biggr\}.
\end{equation*}
\end{corollary}
Application of integral representations (\ref{eq:FStieltjes}), (\ref{eq:pFqgeneral}), (\ref{eq:pFqpsi2-0}), (\ref{eq:pFqsmallp}) (\ref{eq:q+1FqLaplace}), (\ref{eq:qFqLaplace}) and (\ref{eq:q-1Fqcosine}) for investigating the properties of the generalized hypergeometric function $_pF_q$  depends heavily on the positivity of  representing measures, expressed here in terms of Meijer's $G$-function. Sufficient conditions for such positivity are furnished in the next theorem.

\begin{theorem}\label{th:G-positive}
Suppose $A,B\in\R^q$ are such that
\begin{equation}\label{eq:v-def}
v(t)=\sum\limits_{j=1}^{q}(t^{a_j}-t^{b_j})\ge0~~\text{on}~(0,1].
\end{equation}
Then
\begin{equation}\label{eq:G-nonneg}
G^{q,0}_{q,q}\left(t\left|\begin{array}{l}\!\!B\!\!\\\!\!A\!\!\end{array}\right.\right)\ge0~\text{on}~(0,1).
\end{equation}
\end{theorem}

Before giving a proof of this theorem let us remind the reader that a nonnegative function $f$ defined on $(0,\infty)$ is called completely monotone if it has derivatives of all orders and $(-1)^nf^{(n)}(x)\ge0$ for $n\in\N_0$ and $x>0$ \cite[Defintion~1.3]{SSV}. This inequality  is known to be strict unless $f$ is a constant. By the celebrated Bernstein theorem a function is completely monotone if and only if it is the Laplace transform of a nonnegative measure \cite[Theorem~1.4]{SSV}. A positive function $f$ is said to be logarithmically completely monotone if $-(\log{f})'$ is completely monotone \cite[Definition~5.8]{SSV}.  The class of logarithmically completely monotone functions is a proper subset of the class of completely monotone functions.   Their importance stems from the fact that they represent Laplace transforms of infinitely divisible probability distributions, see \cite[Theorem~5.9]{SSV} and  \cite[Section~51]{Sato}.

\textbf{Proof of Theorem~\ref{th:G-positive}.} First note that
$$
G^{q,0}_{q,q}\left(t\left|\begin{array}{l}\!\!B\!\!\\\!\!A\!\!\end{array}\right.\right)=0
$$
for  $t>1$ and all (complex) values of $A$ and $B$. This follows from the fact that for $t>1$  choosing the right loop to be the contour of integration in  (\ref{eq:G-defined}) gives convergent integral according to \cite[Theorem~1.1]{KilSaig}. On the other hand, there are no poles of the integrand inside this contour so that the above equality follows by Cauchy's theorem.
This explains the restriction $t\in(0,1)$ in the formulation of the theorem.  Further, due to the formula
$$
t^{\alpha}G^{q,0}_{q,q}\left(t\left|\begin{array}{l}\!\!B\!\!\\\!\!A\!\!\end{array}\right.\right)
=G^{q,0}_{q,q}\left(t\left|\begin{array}{l}\!\!B+\alpha\!\!\\\!\!A+\alpha\!\!\end{array}\right.\right)
$$
(see \cite[formula~8.2.2.15]{PBM3} or \cite[16.19.2]{NISTCh16}) we can restrict our attention to the case $A,B>0$. Indeed, adding large enough $\alpha$ to $A$ and $B$ neither alters the sign of Meijer's $G$ in (\ref{eq:G-nonneg}) nor the sign of $v(t)$ in (\ref{eq:v-def}).  Adopting the assumption  $A,B>0$ we are in the position to apply \cite[Lemma~2.1]{GrinIsmail} whose particular case (essentially contained already in \cite[Theorem~10]{Alzer}) states that the ratio $x\to\Gamma(A+x)/\Gamma(B+x)$ is logarithmically completely monotone if and only if condition (\ref{eq:v-def}) is satisfied.  Hence, under (\ref{eq:v-def}) this function is also completely monotone. If $\psi>0$ then
$$
\frac{\Gamma(A+x)}{\Gamma(B+x)}=\int_0^{\infty}e^{-xt}G^{q,0}_{q,q}\left(e^{-t}\left|\begin{array}{l}\!\!B\!\!\\\!\!A\!\!\end{array}\right.\right)dt
$$
and the representing measure must be nonnegative by Bernstein's theorem. This measure is unique according to \cite[Proposition~1.2]{SSV}.  Nonnegativity is extended to $\psi=0$ by continuity.  If $\psi<0$ then $v(t)$ cannot be nonnegative on $(0,1]$ since $v(1)=0$ and $v'(1)=-\psi$. $\hfill$$\square$

Condition (\ref{eq:v-def}) is probably also necessary for (\ref{eq:G-nonneg}) at least when $\psi\ge0$. However, this condition is very difficult to verify.  Some sufficient conditions are known for inequality (\ref{eq:v-def}) to hold.  To cite the corresponding results we need to introduce the following terminology \cite[Definition~A.2]{MOA}. It is said that the real vector  $B=(b_1,\ldots,b_q)$ is weakly supermajorized by the real vector $A=(a_1,\ldots,a_q)$ (symbolized as $B\prec^W{A}$) if
\begin{equation}\label{eq:amajorb}
\begin{split}
& 0<a_1\leq{a_2}\leq\cdots\leq{a_q},~~
0<b_1\leq{b_2}\leq\cdots\leq{b_q},
\\
&~\text{and}~\sum\limits_{i=1}^{k}a_i\leq\sum\limits_{i=1}^{k}b_i~~\text{for}~~k=1,2\ldots,q.
\end{split}
\end{equation}
If, in addition, $\psi(=\sum_{i=1}^{q}(b_i-a_i))=0$, than  $B$ is said to be majorized by $A$, or $B\prec{A}$.

It will be convenient to assume that $A$ and $B$ (or $A_i$, $B_i$ when they appear) are ordered ascending whenever they are real.  It follows immediately from a theorem of Tomi\'{c} (see \cite[Proposition~4.B.2]{MOA}) that $v(t)\ge0$ if $B\prec^W{A}$. In the present context this fact was first used by  Alzer \cite[Theorem~10]{Alzer}.  For the particular situation $q=2^n$, $n=0,1,\ldots$, Grinshpan and Ismail \cite[Theorems~1.1,1.2]{GrinIsmail} derived two different sets of sufficient conditions for validity of (\ref{eq:v-def}).  Note that for $p=2$ conditions $\min(a_1,a_2)\leq\min(b_1,b_2)$ and $\psi\ge0$ are necessary and sufficient for validity of (\ref{eq:v-def}). Indeed, for $a_j,b_j>0$ these conditions are clearly equivalent to $B\prec^W{A}$ and sufficiency follows by Tomi\'{c}'s result. Otherwise add positive constant to all of these numbers as in the proof of Theorem~\ref{th:G-positive}. Necessity of the condition $\psi\ge0$ has been also demonstrated in the proof of Theorem~\ref{th:G-positive}.  To prove the necessity of  $\min(a_1,a_2)\leq\min(b_1,b_2)$ assume that $b_1<\min(a_1,a_2,b_2)$
and write $v(t)=t^{b_1}(t^{a_1-b_1}+t^{a_2-b_1}-t^{b_2-b_1}-1)$ which implies that $v(t)<0$ for $t$ near zero.

Combining nonnegativity of $G$-function with representations (\ref{eq:FStieltjes}) and (\ref{eq:qFqLaplace}) we obtain some sufficient conditions for the generalized hypergeometric functions to be completely monotone or logarithmically completely monotone.
\begin{theorem}
Suppose $v(t)\ge0$ on $(0,1]$ and $\sigma>0$.  Then the functions
$$
x\to {_{q+1}F_q}\left(\left.\!\!\begin{array}{c}\sigma,A\\B\end{array}\right|-x\!\right)
~~\text{and}~~
x\to {_{q}F_q}\left(\left.\!\!\begin{array}{c}A\\B\end{array}\right|-x\!\right)
$$
are completely monotone on $(0,\infty)$. In particular, this holds if $B\prec^W{A}$.
\end{theorem}

\begin{theorem}\label{th:qFq-logCM}
Suppose $\sigma>0$ and $v(t)\ge0$ on $(0,1]$ \emph{(}in particular, this holds if $B\prec^W{A}$\emph{)}.  Then the function
$$
x\to x^{-\sigma}{_{q+1}F_q}\left(\left.\!\!\begin{array}{c}\sigma,A\\B\end{array}\right|-\frac{1}{x}\!\right)
$$
is completely monotone on $(0,\infty)$. If $0<\sigma\le1$ then it is logarithmically completely monotone.
\end{theorem}
\textbf{Proof.} By factoring the generalized Stieltjes transform  (\ref{eq:FStieltjes}) into repeated Laplace transforms according to \cite[Theorem~8]{KPJCA} we get
$$
x^{-\sigma}{_{q+1}F_{q}}(\sigma,A;B;-1/x)=\frac{1}{\Gamma(\sigma)}\int\limits_{0}^{\infty}e^{-ux}u^{\sigma-1}\!\!\int\nolimits_{0}^{1}\!e^{-ut}d\rho(t)du
=\frac{1}{\Gamma(\sigma)}\!\int\limits_{0}^{\infty}\!e^{-ux}u^{\sigma-1}{_qF_q}(A;B;-u)du,
$$
where
$$
d\rho(t)=\frac{\Gamma(B)}{\Gamma(A)}G^{q,0}_{q,q}\left(t\left|\begin{array}{l}\!\!B\!\!\\\!\!A\!\!\end{array}\right.\right)\frac{dt}{t}
$$
is nonnegative by Theorem~\ref{th:G-positive} which implies complete monotonicity.  Further, according to \cite[Theorem~51.4]{Sato} a probability distribution is infinitely divisible if it has log-convex density.  The function $u^{\sigma-1}\int_0^1e^{-ut}d\rho(t)$ is log-convex for  $0<\sigma\le1$, since both factors are log-convex (the second factor is log-convex by complete monotonicity).  Thus, the function in the statement of the theorem is the Laplace transform of an infinitely divisible distribution and so is logarithmically completely monotone by \cite[Proposition on p.387]{Alzer} or \cite[Theorem~5.9]{SSV}. $\hfill\square$

By applying the methods of proofs from \cite{KPJMAA,KSJAT} to representations (\ref{eq:pFqgeneral}) and (\ref{eq:pFqpsi2-0}) it is straightforward to get the next two propositions (cf. Theorems~4 and 7 from \cite{KPJMAA}).  The symbol
$A_{1}'$ will denote $A_{1}$ without its maximal element.

\begin{theorem}\label{th:qFq-ratio}
Keep notation and constraints of Theorem~\ref{th:pFqgeneral}  and suppose in addition that $A_1,B_1>0$, $p_2=q_2$ and $\sum_{j=p_1+1}^{p}(t^{a_j}-t^{b_j})\ge0$ \emph{(}or $B_2\prec^W{A_2}$\emph{)}.  Then the function
$$
x\to\frac{{_{p}F_q}\left(\left.\!\!\begin{array}{c}A_1,A_2+\mu\\B_1,B_2+\mu\end{array}\right|-x\!\right)}
{{_{p}F_q}\left(\left.\!\!\begin{array}{c}A_1,A_2\\B_1,B_2\end{array}\right|-x\!\right)}
$$
is monotone decreasing on $(-\infty,0)$ if $p\leq{q}$ or on $(-1,0)$ if $p=q+1$ for every fixed $\mu>0$.
If also $p=q$ and $\sum_{j=1}^{p_1}(t^{a_j}-t^{b_j})\ge0$ \emph{(}or $B_1\prec^W{A_1}$\emph{)}, then the above quotient decreases on the whole real line.  If $p=q+1$ and $\sum_{j=1}^{q_1}(t^{a_j}-t^{b_j})\ge0$ \emph{(}or $B_1\prec^W{A_{1}'}$\emph{)}, then the above quotient decreases on $(-1,\infty)$.
\end{theorem}

\begin{theorem}\label{th:qFq-logconv}
Keep notation and constraints of Theorem~\ref{th:pFqgeneral} and suppose in addition that $A_1,B_1>0$, $p_2=q_2$ and $\sum_{j=p_1+1}^{p}(t^{a_j}-t^{b_j})\ge0$ \emph{(}or $B_2\prec^W{A_2}$\emph{)}.  Then the function
$$
\mu\to{_{p}F_q}\left(\left.\!\!\begin{array}{c}A_1,A_2+\mu\\B_1,B_2+\mu\end{array}\right|-x\!\right)
$$
is log-convex on $(0,\infty)$ for each fixed $x\in(-\infty,0)$ if $p\leq{q}$ or $x\in(-1,0)$ if
 $p=q+1$. If also $p=q$ and $\sum_{j=1}^{p_1}(t^{a_j}-t^{b_j})\ge0$ \emph{(}or $B_1\prec^W{A_1}$\emph{)}, then log-convexity holds for each real $x$, while for $p=q+1$  and $\sum_{j=1}^{q_1}(t^{a_j}-t^{b_j})\ge0$ \emph{(}or $B_1\prec^W{A_{1}'}$\emph{)} log-convexity holds  for each fixed  $x\in(-1,\infty)$.
\end{theorem}

\textbf{Remark.}  It is easy to see that conditions  $B_1\prec^W{A_1}$ and  $B_2\prec^W{A_2}$ imply $B\prec^W{A}$
(for these relations to make sense one has to assume that  $p_1=q_1$ and $p_2=q_2$).  For this reason Theorems~\ref{th:qFq-ratio} and \ref{th:qFq-logconv} are the strongest in some informal sense when $p_1=q_1=0$, i.e. for the functions
$$
x\to{_{q}F_q}\left(\left.\!\!\begin{array}{c}A+\mu\\B+\mu\end{array}\right|-x\!\right)\!\Big/\!
{_{q}F_q}\left(\left.\!\!\begin{array}{c}A\\B\end{array}\right|-x\!\right)~~\text{and}~~\mu\to {_{q}F_q}\left(\left.\!\!\begin{array}{c}A+\mu\\
B+\mu\end{array}\right|-x\!\right).
$$

\section{Inequalities for the Kummer and Gauss type functions}

In Theorem~16 of his paper \cite{Luke} Luke gave two-sided bounds for the function ${_qF_q}(A;B;x)$ under the restrictions
$b_i\ge{a_i}>0$, $i=1,2,\ldots,q$. He indicated that these bounds are ''easily proved'' without providing such proofs.  In this section we will supply  two different proofs of Luke's inequalities, one valid for positive values of the argument $x$ and the other valid for all real $x$. In this way we substantially relax Luke's conditions.  For negative argument values our conditions are given in terms of nonnegativity of $v(t)$ or weak majorization $B\prec^W{A}$. For positive argument values the conditions can be weakened further and are given in terms of elementary symmetric polynomials, defined by
$$
e_k(x_1,\ldots,x_q)=\sum\limits_{1\leq{j_1}<{j_2}\cdots<{j_k}\leq{q}}x_{j_1}x_{j_2}\cdots{x_{j_k}},~~k=1,2,\ldots,q.
$$

\begin{theorem}\label{th:qFqpositive-x}
Suppose
\begin{equation}\label{eq:symmetricba}
\frac{e_q(b_1,\ldots,b_q)}{e_q(a_1,\ldots,a_q)}\geq
\frac{e_{q-1}(b_1,\ldots,b_q)}{e_{q-1}(a_1,\ldots,a_q)}\geq\cdots\geq
\frac{e_1(b_1,\ldots,b_q)}{e_1(a_1,\ldots,a_q)}\geq{1}
\end{equation}
and each elementary symmetric polynomial above is nonnegative.  Then
\begin{equation}\label{eq:qFq-Luke}
e^{f_1x}\leq{_qF_q}(A;B;x)\leq1-f_1+f_1e^x~~\text{for}~~x\ge0,
\end{equation}
where $f_1=\prod_{i=1}^{q}(a_i/b_i)$. Moreover, the upper bound holds true if each fraction in \emph{(\ref{eq:symmetricba})} is merely not less than $1$.
\end{theorem}
\textbf{Remark.}  Note that conditions (\ref{eq:symmetricba}) are strictly weaker than $B\prec^W{A}$,
as we demonstrated in \cite[Lemma~3]{KPJMAA}.

\noindent\textbf{Proof.} Denote by $f_n=\prod_{i=1}^{q}[(a_i)_n/(b_i)_n]$ the coefficient at $x^n/n!$ in power series expansion (\ref{eq:pFqdefined}) of ${_qF_q}(A;B;x)$.   Then conditions
$$
e_{i}(b_1,\ldots,b_q)\ge e_{i}(a_1,\ldots,a_q),~~~i=1,\ldots,q,
$$
(i.e. each fraction in (\ref{eq:symmetricba}) is not less than $1$) imply that
$$
\frac{f_{n+1}}{f_{n}}=R(n)=\prod_{i=1}^{q}\frac{a_i+n}{b_i+n}\le1,
$$
since $e_{q-i}(a_1,\ldots,a_q)$ ($e_{q-i}(b_1,\ldots,b_q)$) is the coefficient of $n^i$ in the polynomial in the numerator (denominator) of $R(n)$.  Thus $f_{n+1}\le{f_{n}}$, so that  $f_n\le{f_1}$ for $n=1,2,\ldots$. Consequently,
for $x\ge0$ we get
$$
{_qF_q}(A;B;x)=1+\sum\limits_{n=1}^{\infty}f_n\frac{x^n}{n!}
=1+f_1\sum\limits_{n=1}^{\infty}\frac{f_n}{f_1}\frac{x^n}{n!}\le1+f_1\sum\limits_{n=1}^{\infty}\frac{x^n}{n!}=1-f_1+f_1e^x.
$$
Further, under conditions (\ref{eq:symmetricba}) the function $R(x)$ defined above is increasing according to \cite[Lemma~2]{KSJMAA}. This leads to the following inequalities ($k\ge0$):
$$
R(0)=\prod_{i=1}^{q}\frac{a_i}{b_i}\leq\prod_{i=1}^{q}\frac{a_i+k}{b_i+k}=R(k)~~\Rightarrow~~
(f_1)^n=\prod_{i=1}^{q}\frac{(a_i)^n}{(b_i)^n}\leq\prod_{i=1}^{q}\frac{(a_i)_n}{(b_i)_n}=f_n,~n=1,2,\ldots
$$
Consequently,
$$
{_qF_q}(A;B;x)=1+\sum\limits_{n=1}^{\infty}f_n\frac{x^n}{n!}\ge 1+\sum\limits_{n=1}^{\infty}(f_1)^n\frac{x^n}{n!}=e^{f_1x},
$$
which completes the proof.\hfill $\square$

\noindent\textbf{Remark}. Inequalities (\ref{eq:qFq-Luke}) can be refined to the estimates
\begin{equation}\label{eq:qFqLuke-refined}
1+\frac{f_1^2}{f_2}(e^{(f_2/f_1)x}-1)\leq{_qF_q}(A;B;x)\leq1-f_2+(f_1-f_2)x+f_2e^x
\end{equation}
valid for $x\ge0$ under conditions of Theorem~\ref{th:qFqpositive-x}.  Indeed, the upper bound is obtained by writing
$$
{_qF_q}(A;B;x)=1+f_1x+f_2\sum\limits_{n=2}^{\infty}\frac{f_n}{f_2}\frac{x^n}{n!}\le
1+f_1x+f_2\sum\limits_{n=2}^{\infty}\frac{x^n}{n!}=1-f_2+(f_1-f_2)x+f_2e^x,
$$
where we used $f_{n+1}\le f_{n}$ for $n=2,3,\ldots$ provided that each fraction in (\ref{eq:symmetricba}) is not less than $1$.
To prove  the lower bound we note that under conditions (\ref{eq:symmetricba}) we have $(f_2/f_1)^{n-1}\leq f_n/f_1$ for $n=2,3,\ldots$ by the increase of $R(x)$.  Then
$$
{_qF_q}(A;B;x)=1+f_1x+f_1\sum\limits_{n=2}^{\infty}\frac{f_n}{f_1}\frac{x^n}{n!}\ge
1+f_1x+f_1\sum\limits_{n=2}^{\infty}\left(\frac{f_2}{f_1}\right)^{n-1}\frac{x^n}{n!}=1+\frac{f_1^2}{f_2}(e^{(f_2/f_1)x}-1).
$$
Similar trick can be applied to separate further terms.

\begin{corollary} Suppose $\sigma>0$ and hypotheses of Theorem~\ref{th:qFqpositive-x} are satisfied. Then
for $0\le{x}<1$
$$
\frac{1}{(1-f_1x)^{\sigma}}\leq{_{q+1}F_q}(\sigma,A;B;x)\leq1-f_1+\frac{f_1}{(1-x)^{\sigma}}
$$
and
$$
1-\frac{f_1^2}{f_2}+\frac{f_1^2}{f_2(1-f_2x/f_1)^{\sigma}}\leq{_{q+1}F_q}(\sigma,A;B;x)\leq
1-f_2+\sigma(f_1-f_2)x+\frac{f_2}{(1-x)^{\sigma}}.
$$
\end{corollary}
\textbf{Proof.}  Following Luke's idea from \cite{Luke}, write the bounds (\ref{eq:qFq-Luke}) for ${_{q}F_q}(A;B;t)$, multiply by $e^{-ty}t^{\sigma-1}$ and integrate using
$$
\int\nolimits_{0}^{\infty}e^{-ty}t^{\sigma-1}{_{q}F_q}(A;B;t)dt=y^{-\sigma}\Gamma(\sigma){_{q+1}F_q}(\sigma,A;B;1/y).
$$
It remains to write $x=1/y$ in the resulting inequality and simplify to get the first inequality.  The second inequality is proved by applying the same trick to (\ref{eq:qFqLuke-refined}).$\hfill\square$

\begin{theorem}\label{th:qFq-Jensen}
Suppose $A,B>0$ and $\sum_{j=1}^{q}(t^{a_j}-t^{b_j})\ge0$ \emph{(}or $B\prec^W A$\emph{)}. Then
$$
e^{-f_1x}\le{_qF_q}(A;B;-x)\le1-f_1+f_1e^{-x}
$$
for all real $x$.
\end{theorem}
\textbf{Proof.}  According to the integral form of  Jensen's inequality  \cite[Chapter~I, formula (7.15)]{MPF}
\begin{equation}\label{eq:Jensen}
\varphi\left(\int\limits_a^{b}f(s)d\mu(s)\Bigg/\int\limits_a^{b}d\mu(s)\right)
\leq\int\limits_a^{b}\varphi(f(s))d\mu(s)\Bigg/\int\limits_a^{b}d\mu(s)
\end{equation}
if  $\varphi$ is convex and $f$ is integrable with respect to a nonnegative measure $\mu$.
Put $\varphi_x(y)=e^{-xy}$, $f(s)=s$ and
$$
d\mu(s)=\frac{\Gamma(B)}{\Gamma(A)}G^{q,0}_{q,q}\left(s\left|\begin{array}{l}\!\!B\!\!\\\!\!A\!\!\end{array}\right.\right)\frac{ds}{s}.
$$
Then
$$
\int\limits_0^{1}d\mu(s)=1,~~~~~\int\limits_0^{1}f(s)d\mu(s)=\prod\limits_{i=1}^{q}\frac{a_i}{b_i}=f_1,~~~\int\limits_0^{1}\varphi_x(f(s))d\mu(s)={_qF_q}(A;B;-x).
$$
The last equality is a rewriting of (\ref{eq:qFqLaplace}). This proves the lower bound. To demonstrate the upper bound we will apply the converse Jensen inequality due to Lah and Ribari\'{c}, which reads as follows. Set
$$
A(g)=\left.\int\limits_m^{M}g(s)d\tau(s)\middle/ \int\limits_m^{M}d\tau(s)\right.,
$$
where $\tau$ is a nonnegative measure and  $g$ is a continuous function.  If $-\infty<m<M<\infty$ and $\varphi$ is convex on $[m,M]$ then according to \cite[Theorem~3.37]{PPT}
$$
(M-m)A(\varphi(g))\leq (M-A(g))\varphi(m)+(A(g)-m)\varphi(M).
$$
Setting $\varphi_x(t)=e^{-xt}$, $d\tau(s)=d\mu(s)$,  $g(s)=s$ and $[m,M]=[0,1]$, we arrive at the upper bound of the theorem.\hfill$\square$

\begin{corollary} Suppose $\sigma>0$ and hypotheses of Theorem~\ref{th:qFq-Jensen} are satisfied. Then
$$
\frac{1}{(1+f_1x)^{\sigma}}\leq{_{q+1}F_q}(\sigma,A;B;-x)\leq1-f_1+\frac{f_1}{(1+x)^{\sigma}}
$$
for $x\ge0$.
\end{corollary}

\noindent\textbf{Proof.} Multiply inequality (\ref{eq:qFq-Luke}) written for ${_qF_q}(A;B;-xt)$ by $e^{-t}t^{\sigma-1}$ and integrate using the formula
$$
\int_0^{\infty}e^{-t}t^{\sigma-1}{_qF_q}(A;B;-xt)dt=\Gamma(\sigma){_{q+1}F_q}(\sigma, A;B;-x).~~~~~~\square
$$

\section{Inequalities for the Bessel type functions}

First, we will find an upper bound in the general situation $p<q$.  As before the symbol $f_n$ will denote the coefficient at $x^n/n!$ in the series representation
(\ref{eq:pFqdefined}), i.e.
$$
f_n=\frac{\prod_{i=1}^{p}(a_i)_n}{\prod_{i=1}^{q}(b_i)_n}=\frac{(A)_n}{(B)_n}~\text{for}~n=0,1,\ldots
$$
\begin{theorem}\label{th:upperpFq}
Suppose $p<q$.  If
\begin{equation}\label{eq:ebea}
e_{q-i}(b_1,\ldots,b_q)\ge e_{p-i}(a_1,\ldots,a_p),~~~i=0,1,\ldots,p,
\end{equation}
then for $x\ge0$
\begin{equation*}
{_pF_q}(A;B;x)\leq1-f_1+f_1e^x.
\end{equation*}
If
\begin{equation}\label{eq:R-decreasing}
\frac{e_q(b_1,\ldots,b_q)}{e_{p}(a_1,\ldots,a_{p})}\leq
\frac{e_{q-1}(b_1,\ldots,b_q)}{e_{p-1}(a_1,\ldots,a_{p})}\leq\cdots
\leq\frac{e_{q-p+1}(b_1,\ldots,b_q)}{e_{1}(a_1,\ldots,a_{p})}\leq
e_{q-p}(b_1,\ldots,b_q),
\end{equation}
then for $x\ge0$
$$
{_pF_q}(A;B;x)\le e^{f_1x}.
$$
\end{theorem}
\textbf{Proof.}  The proof of the first upper bound repeats the proof of the upper bound (\ref{eq:qFq-Luke})
in Theorem~\ref{th:qFqpositive-x}.  To demonstrate the second bound note that for $p<q$ the function
$$
R(x)=\frac{\prod_{i=1}^{p}(a_i+x)}{\prod_{i=1}^{q}(b_i+x)}
$$
is decreasing under conditions (\ref{eq:R-decreasing}) according to \cite[p.394]{KSJMAA}
which implies that
$$
f_n=\frac{\prod_{i=1}^{p}(a_i)_n}{\prod_{i=1}^{q}(b_i)_n}\leq\frac{\prod_{i=1}^{p}(a_i)^n}{\prod_{i=1}^{q}(b_i)^n}=(f_1)^n.
$$
Hence,
$$
{_pF_q}(A;B;x)=1+\sum\limits_{n=1}^{\infty}f_n\frac{x^n}{n!}\le 1+\sum\limits_{n=1}^{\infty}(f_1)^n\frac{x^n}{n!}=e^{f_1x}.~~~~~~~~~~~~~~~\square
$$

According to the asymptotic formula \cite[16.11.8]{NISTCh16},
$$
{_{q-1}F_{q}}\left(\left.\!\!\begin{array}{l}A\\B\end{array}\right|x\right)=
\frac{\Gamma(b_1)\cdots\Gamma(b_q)}{2\sqrt{\pi}\Gamma(a_1)\cdots\Gamma(a_{q-1})}x^{\nu}e^{2\sqrt{x}}\left(1+\frac{d_1}{\sqrt{x}}+O(x^{-1})\right)~\text{as}~x\to+\infty,
$$
where $\nu=\frac{1}{2}\sum_{i=1}^{q-1}a_i-\frac{1}{2}\sum_{i=1}^{q}b_i+1/4$. Hence,  the upper bounds of Theorem~\ref{th:upperpFq} are very wrong in asymptotic order.   In the most important case $p=q-1$ we can do much better.

\begin{theorem}
Suppose $A,B>0$ \emph{(}understood element-wise\emph{)}. Then for $x\ge0$
\begin{equation}\label{eq:pFqlower}
e^{\sqrt{4x+c^2}-c}\left(\frac{1}{2}+\frac{1}{2c}\sqrt{4x+c^2}\right)^{-c}\le{_{q-1}F_q}(A;B;x),
\end{equation}
where $c>0$ is given by
\begin{equation}\label{eq:c-max}
c=\max\limits_{i\in\{1,2,\ldots,q\}}\left[\frac{e_i(b_1,b_2,\ldots,b_q)-e_i(a_1,a_2,\ldots,a_p)}{e_{i-1}(a_1,a_2,\ldots,a_p)}\right],
\end{equation}
$p=q-1$ and $e_q(a_1,a_2,\ldots,a_p)=0$.
\end{theorem}
\textbf{Proof.} Suppose we could find a number $c$ such that
\begin{equation}\label{eq:fn-cn}
f_n=\frac{(a_1)_n\cdots(a_p)_n}{(b_1)_n\cdots(b_q)_n}\ge\frac{1}{(c)_n}~~\text{for}~~n=1,2,\ldots
\end{equation}
Then for $x\ge0$ ($p=q-1$):
\begin{equation}\label{eq:q-1Fqlower}
{_pF_q}(A;B;x)=1+\sum\limits_{n=1}^{\infty}\frac{(a_1)_n\cdots(a_p)_n}{(b_1)_n\cdots(b_q)_n}\frac{x^n}{n!}\ge 1+\sum\limits_{n=1}^{\infty}\frac{1}{(c)_n}\frac{x^n}{n!}={_0F_1}(-;c;x).
\end{equation}
Further, we can use some known  lower bounds for the function ${_0F_1}(-;c;x)$ (which is equal to the modified Bessel function $I_{c-1}$ up to a simple multiplier) to derive lower bounds for ${_pF_q}(A;B;x)$ in terms of elementary functions.
For (\ref{eq:fn-cn}) to hold it suffices to satisfy $f_1c\ge1$ and
$$
\frac{f_{n+1}(c)_{n+1}}{f_n(c)_{n}}=\frac{(a_1+n)\cdots(a_p+n)(c+n)}{(b_1+n)\cdots(b_q+n)}\ge1,~~n=1,2,\ldots
$$
In turn, the above inequality  holds if  (recall that $q=p+1$)
$$
e_i(a_1,a_2,\ldots,a_p,c)\ge e_i(b_1,b_2,\ldots,b_q),~~~~i=1,2,\ldots,q,
$$
or
$$
e_i(a_1,a_2,\ldots,a_p)+ce_{i-1}(a_1,a_2,\ldots,a_p) \ge e_i(b_1,b_2,\ldots,b_q),~~~~i=1,2,\ldots,q.
$$
To satisfy these $q$ inequalities we need to choose $c$ by
$$
c=\max\limits_{i\in\{1,2,\ldots,q\}}\left[\frac{e_i(b_1,b_2,\ldots,b_q)-e_i(a_1,a_2,\ldots,a_p)}{e_{i-1}(a_1,a_2,\ldots,a_p)}\right].
$$
Here $e_0=1$ and $e_q(a_1,a_2,\ldots,a_p)=0$. Due to the last identity we get $c>0$  for any positive arrays $A$ and $B$.
Hence, the problem reduces to finding good bounds for ${_0F_1}(-;c;x)$ for $x,c>0$.  Numerically best bounds are contained in \cite[formula (11)]{Amos} in terms of the ratio $I_{\nu+1}/I_{\nu}$ of the modified Bessel functions
$$
I_{\nu}(x)=(x/2)^{\nu}[\Gamma(\nu+1)]^{-1}{_0F_1}(-;\nu+1;x^2/4).
$$
When rewritten in terms of the logarithmic derivative of ${_0F_1}(-;c;x)$ these bounds read
$$
\frac{2}{\sqrt{4x+c^2}+c}\le\frac{{_0F_1}\!'(-;c;x)}{{_0F_1}(-;c;x)}=\frac{{_0F_1}(-;c+1;x)/c}{{_0F_1}(-;c;x)}\le\frac{2}{\sqrt{4x+(c+1)^2}+c-1},
$$
where the derivative formula ${_0F_1}\!'(-;c;x)={_0F_1}(-;c+1;x)/c$ has been used.
Employing the evaluation
$$
\int\limits_{0}^{x}\frac{2dt}{a+\sqrt{4tq+b^2}}=\frac{1}{q}\sqrt{4xq+b^2}-\frac{a}{q}\ln\frac{a+\sqrt{4xq+b^2}}{a+b}-\frac{b}{q},
$$
we can integrate the above inequalities to obtain
\begin{multline*}
\sqrt{4x+c^2}-c\log\frac{c+\sqrt{4x+c^2}}{2c}-c\le\log({_0F_1}(-;c;x))
\\
\le\sqrt{4x+(c+1)^2}-(c-1)\log\frac{c-1+\sqrt{4x+(c+1)^2}}{2c}-(c+1).
\end{multline*}
Taking exponentials yields:
\begin{equation}\label{eq:0F1bounds}
e^{\sqrt{4x+c^2}-c}\!\left(\frac{1}{2}+\frac{1}{2c}\sqrt{4x+c^2}\right)^{\!\!-c}\!\!\!\!\le{_0F_1}(-;c;x)\le e^{\sqrt{4x+(c+1)^2}-c-1}\!\left(\frac{c-1}{2c}+\frac{1}{2c}\sqrt{4x+(c+1)^2}\right)^{\!\!1-c}.
\end{equation}
Combining the lower bound here with (\ref{eq:q-1Fqlower}) proves the theorem. $\hfill\square$
\begin{theorem}
Suppose $A,B>0$ \emph{(}understood element-wise\emph{)} and $d$ given by
\begin{equation}\label{eq:d-min}
d=\min\limits_{i\in\{1,2,\ldots,q\}}\left[\frac{e_i(b_1,b_2,\ldots,b_q)-e_i(a_1,a_2,\ldots,a_p)}{e_{i-1}(a_1,a_2,\ldots,a_p)}\right]
\end{equation}
is positive. Here $p=q-1$, $e_0=1$ and $e_q(a_1,a_2,\ldots,a_p)=0$.  Then for $x\ge0$
\begin{equation}\label{eq:pFqupper}
{_{q-1}F_q}(A;B;x)\le e^{\sqrt{4x+(d+1)^2}-d-1}\left(\frac{d-1}{2d}+\frac{1}{2d}\sqrt{4x+(d+1)^2}\right)^{1-d}.
\end{equation}
\end{theorem}
\textbf{Proof.}  If we could find $d$ such that
$$
f_n(d)_n=\frac{(a_1)_n\cdots(a_p)_n(d)_n}{(b_1)_n\cdots(b_q)_n}\le1~~\text{for}~~n=1,2,\ldots,
$$
then
$$
{_pF_q}(A;B;x)=1+\sum\limits_{n=1}^{\infty}\frac{(a_1)_n\cdots(a_p)_n}{(b_1)_n\cdots(b_q)_n}\frac{x^n}{n!}\le 1+\sum\limits_{n=1}^{\infty}\frac{1}{(d)_n}\frac{x^n}{n!}={_0F_1}(-;d;x).
$$
Application of the upper bound from (\ref{eq:0F1bounds}) to the above inequality will prove (\ref{eq:pFqupper}).
To find such $d$ it suffices to satisfy $f_1d\le1$ and
$$
\frac{f_{n+1}(d)_{n+1}}{f_n(d)_{n}}=\frac{(a_1+n)\cdots(a_p+n)(d+n)}{(b_1+n)\cdots(b_q+n)}\le1.
$$
In turn, the above inequality  holds if (recall that $q=p+1$)
$$
e_i(a_1,a_2,\ldots,a_p,d)\le e_i(b_1,b_2,\ldots,b_q),~~~~i=1,2,\ldots,q,
$$
or
$$
e_i(a_1,a_2,\ldots,a_p)+de_{i-1}(a_1,a_2,\ldots,a_p) \le e_i(b_1,b_2,\ldots,b_q),~~~~i=1,2,\ldots,q.
$$
To satisfy these $q$ inequalities we need to choose $d$ by (\ref{eq:d-min}).$\hfill\square$

\bigskip

\textbf{Acknowledgements.} This work has been supported by the Russian Science Foundation under project 14-11-00022.

\end{document}